\documentclass[a4paper,12pt]{article}

\usepackage{tgtermes}

\usepackage[utf8]{inputenc}
\usepackage{amsmath}
\usepackage{amsfonts}
\usepackage{amssymb}
\usepackage{fontenc}
\usepackage{graphicx}

\usepackage[dvips]{hyperref}
\usepackage{xcolor}

\usepackage{mathrsfs}

\providecommand{\keywords}[1]{\textbf{Keywords }#1}

\def\upd{\,{\rm d}}
\newcommand{\be}{\begin{eqnarray}}
\newcommand{\ee}{\end{eqnarray}}
\newcommand{\nnee}{\nonumber\ee}
\newcommand{\Tr}{\,{\rm Tr}\,}
\newcommand{\dom}{\,{\rm dom}\,}

\def\beginproof{\par\strut\vskip 0.0cm\noindent{\bf Proof}\par}
\def\endproof{\par\strut\hfill$\square$\par\vskip 0.2cm}

\def\Co{{\mathbb C}}

\def\Io{{\mathbb I}}

\def\Ro{{\mathbb R}}

\def\Vo{{\mathbb V}}

\newtheorem{theorem}{Theorem}[section]

\newtheorem{proposition}[theorem]{Proposition}
\newtheorem{definition}[theorem]{Definition}

\newtheorem{corollary}[theorem]{Corollary}

\title{Exponential arcs in the manifold of\\
vector states on a $\sigma$-finite\\
von Neumann algebra
}

\author{ Jan Naudts \\
Physics Department, Universiteit Antwerpen,\\
Universiteitsplein 1, 2610 Antwerpen Belgium\\
e-mail: jan.naudts@uantwerpen.be \\
orcid.org: 0000-0002-4646-1190
}

\date{}

\begin{document}
\maketitle

\begin{abstract}
This paper introduces the notion of 
exponential arcs in Hilbert space and of
exponential arcs connecting vector states
on a sigma-finite von Neumann algebra in its standard representation.
Results from Tomita-Takesaki theory form an essential ingredient.
Starting point is a non-commutative Radon-Nikodym theorem that
involves positive operators affiliated with the commutant algebra.
It is shown that exponential arcs are differentiable and that
parts of an exponential arc are again exponential arcs.
Special cases of probability theory and of quantum probability are used
to illustrate the approach.
\end{abstract}

\keywords{
Exponential arc, Exponential family, Tomita-Takesaki theory, Information Geometry,
Probability theory, Quantum Probability}

\section{Introduction}

Araki \cite{AH74} introduces $\alpha$-families of cones in Hilbert space $\mathscr H$.
They are characterized by a $\sigma$-finite von Neumann algebra $\cal A$ of operators in $\mathscr H$
and a vector $x$ in $\mathscr H$ which is cyclic and separating for $\cal A$.
The parameter $\alpha$ satisfies $0\le\alpha\le 1/2$. The $\alpha$-cone and the $(\frac 12-\alpha)$-cone
are each others dual. Hence, the cone with $\alpha=1/4$ is self-dual. It is called the
natural positive cone and it plays a special role
in Tomita-Takesaki theory \cite{TM70}.
For each value of $\alpha$ one can define generalized Radon-Nikodym operators.
In the case $\alpha=0$ they belong to the algebra $\cal A$, in the case $\alpha=1/2$ they belong
to the commutant algebra ${\cal A}'$. The present paper is restricted to the $\alpha=1/2$-case
in the belief that $\alpha=1/2$ simplifies the introduction in a non-commutative context
of the notion of models belonging to an exponential family.

Amari \cite{AS85,AN00,AS16} introduces an $\alpha$-collection of geometries on manifolds of non-degenerate
probability distributions. The $\alpha$ and $-\alpha$ geometries are each other dual. The self-dual
case $\alpha=0$ involves the Levi-Civita connection. Models belonging to an exponential
family have a flat geometry in the $\alpha=1$ connection.

The analogy between Amari's $\alpha$-geometries and Araki's $\alpha$-cones is
striking but may be superficial. Both theories have in common that the introduction
of the $\alpha$-dependence enriches the theory. In addition, the self-dual choice of $\alpha$ is
technically more demanding. The study of statistical models belonging to an exponential
family is more easy when using a flat geometry, this is, $\alpha=1$.
Goal of the present paper is
to explore the use of Araki's $\alpha$-cones with $\alpha=1/2$
in the context of ongoing efforts to generalize Information Geometry to a non-commutative setting. 

A main concept in Information Geometry \cite{AN00,AS16,AJLS18}
is the notion of a model belonging to an exponential family.
In Statistical Physics such a model is called a Boltzmann-Gibbs distribution.
The generalization \cite{HH93,HH97,HP97,PS96,OSA96,AN00b} to quantum models 
is feasible, at least when the Hilbert space of the model is finite-dimensional
and the states on the algebra can be described by density matrices.
As soon as it is infinite dimensional technical problems appear on top of the 
problems due to non-commutativity. 
Operators on Hilbert space need not anymore to be continuous and their domain of definition
is only a subspace. Selfadjoint operators play a crucial role. Non-uniqueness of
selfadjoint extensions is a serious problem.
Early attempts to tackle these problems include \cite{GS97,SRF04a,SRF04b,JA01,JA06}.

Cones of positive-definite matrices are studied in \cite{OA04,UO04,OA19,KU19}.
These cones are self-dual for the Hilbert-Schmidt inner product.
With each such cone corresponds a Euclidean Jordan algebra. The product of this algebra
is commutative. This simplifies the introduction of Amari's dually flat geometry.
It is not immediately clear how to generalize this approach to infinite dimensions.

The recent approach of \cite{CIJM19,CIJM19a} relies on Lie theory for the
group of bounded operators with bounded inverse and does not yet attack the
above mentioned problems with unbounded operators.

The present approach involving the commutant algebra ${\cal A}'$
has already been explored \cite{NJ18,NJ19b} in the familiar context of density matrices.
It is also studied in \cite{NJ19} where the exponential function is replaced by a
function with linear growth at $+\infty$. The latter idea is
tried out in the context of commutative probability theory with the intent to
better control problems of unboundedness \cite{NN12,NN19,MP17,MP19}.
It is not used in the present paper
because Tomita-Takesaki theory relies heavily on the properties of the exponential function.

\paragraph{Structure of the paper}
The next section reviews elements of Tomita-Takesaki-theory and fixes some of the notations.
Section \ref {sect:Hsp} introduces exponential arcs in Hilbert space.
Section \ref {sect:ssp} introduces exponential arcs in the space of vector states.
Section \ref {sect:special} considers special cases.
The final section contains a short discussion.

\section{Known and less-known results}

Throughout the text $\cal A$ is a fixed $\sigma$-finite von Neumann algebra of operators on a Hilbert 
space $\mathscr H$ and $x$ is a fixed normalized vector in $\mathscr H$ that is 
cyclic and separating for $\cal A$. 
Each vector $y$ in $\mathscr H$ defines a linear functional $\omega_y$ on $\cal A$ by
\be
\omega_y(A)&=&(Ay,y),
\qquad A\in{\cal A}
\label{manif:xidef}
\ee
If $||y||=1$ then $\omega_y$ is a {\em vector state} on $\cal A$. 
The set of vector states of $\cal A$ is denoted ${\Vo}$.
A state $\omega$ on $\cal A$ is said to be {\em faithful} if $\omega(A^*A)=0$
implies $A=0$. The vector state $\omega_x$ is faithful if $x$ is separating for $\cal A$.
The subset of faithful vector states is denoted $\mathring\Vo$.

\subsection{Cones in Hilbert space}

The closure of the map
\be
Ax\mapsto A^*x,
\qquad A\in{\cal A}
\nnee
is denoted $S_x$. 
Let $S_x=J_x|\Delta_x|^{1/2}$ be its polar decomposition.
The {\em modular automorphisms} $\tau_t$, $t\in\Ro$, are defined by
\be
\tau_t(A)=\Delta_x^{it}A\Delta_x^{-it},
\qquad A\in{\cal A}.
\nnee
The anti-linear isometry $J_x$ satisfies $J_x^*=J_x$ and $J_x^2=\Io$.

A cone ${\mathscr P}_x$ associated with the couple ${\cal A},x$ is defined 
by taking the closure of the set of vectors
\be
\{AJ_xAx:\,A\in{\cal A}\}.
\nnee
If $x$ is cyclic and separating for $\cal A$ then 
${\mathscr P}_x$ is independent \cite{AH74} of the choice of $x$.
It is called the {\em natural positive cone}.
It has the properties that 
\be 
{\mathscr P}_x={\mathscr P}_y:= {\mathscr P}
\mbox{ and } 
J_x=J_y:=J 
\quad\mbox{ for all } x,y\in {\mathscr P}.
\nnee
For any $x$ in ${\mathscr P}$ and $A$ in $\cal A$ is $(AJAx,x)\ge0$ with equality if and only if
$A=0$.

Introduce the notations
\be
{\mathscr C}_x&=&\{y\in{\mathscr H}:\, (y,Ax)\ge 0 \mbox{ holds for all }A\ge 0, A\in{\cal A}\}
\nnee
and
\be
{\mathscr C}^*_x&=&\{y\in{\mathscr H}:\, (y,A'x)\ge 0 \mbox{ holds for all }A'\ge 0, A'\in{\cal A'}\}.
\nnee
The sets ${\mathscr C}_x$ and ${\mathscr C}^*_x$ are cones
closed in the norm topology of $\mathscr H$.
A vector $y$ belongs to ${\mathscr C}_x$ if and only if $J_xy$ belongs to ${\mathscr C}^*_x$.
The cones ${\mathscr C}_x$ and ${\mathscr C}^*_x$ form a dual pair
in the sense that
\be
{\mathscr C}_x&=&\{y\in{\mathscr H}:\, (y,z)\ge 0 \mbox{ for all }z\in {\mathscr C}^*_x\}
\nnee
and
\be
{\mathscr C}^*_x&=&\{y\in{\mathscr H}:\, (y,z)\ge 0 \mbox{ for all }z\in {\mathscr C}_x\}.
\nnee

The notations ${\mathscr C}_x$ and ${\mathscr C}^*_x$ are used rather than ${\mathscr P}^\flat_x$,
respectively ${\mathscr P}^\sharp_x$, because further indices will be added later on.
In \cite{AH74} the notation $V^\alpha$ with $\alpha=1/2$, respectively $\alpha=0$ is used. 
The cone $V^{1/4}$ is then the natural positive cone.

The next result is important for what follows.

\begin{proposition}
\label{prop:coneposop}
A vector $y\in{\mathscr H}$ belongs to ${\mathscr C}_x$ if and only if there exists
a non-negative selfadjoint operator $Y$
affiliated with ${\cal A}'$ such that $y=Yx$.
\end{proposition}

For a proof, see for instance Proposition 2.5.27 of \cite{BR79}.

\begin{corollary}
\label{poscon:corol1}
Any vector $y$ of ${\mathscr C}_x$ satisfies $S^*_xy=y$.
In particular, one has 
\be
(y,Ax)&=&(x,Ay)
\quad\mbox{for all}\quad
A\in {\cal A}.
\nnee
\end{corollary}

\beginproof
One has for all $A$ in ${\cal A}$ that
\be
(S_xAx,y)=(A^*x,Yx)=(Yx,Ax)=(y,Ax).
\nnee
Hence, $y$ is in the domain of $S^*_x$ and satisfies $S^*_xy=y$.
One has
\be
(y,Ax)=(y,S_xA^*x)=(A^*x,S^*_xy)=(A^*x,y)=(x,Ay).
\nnee

\endproof

\begin{definition}
Let ${\mathscr C}^1_x$ denote the subset of normalized elements of ${\mathscr C}_x$.
Let $\breve{\mathscr C}_x={\cal A}'x$ and let $\breve{\mathscr C}^1_x$ denote the subset
 of normalized vectors in $\breve{\mathscr C}_x$.
\end{definition}

\begin{proposition}
\label{cones:prop:dense}
The set $\breve{\mathscr C}_x$ is dense in ${\mathscr C}_x$;
The set $\breve{\mathscr C}^1_x$ is dense in ${\mathscr C}^1_x$. 
\end{proposition}

\beginproof
Let $y\in {\mathscr C}_x$.
By Proposition \ref {prop:coneposop} there exists a non-negative operator $Y$ affiliated with the
commutant ${\cal A}'$ such that $y=Yx$.
The operator $Y$ can be approximated by non-negative operators
$Y_n$ belonging to the commutant ${\cal A}'$ in such a way that the vectors $y_n=Y_nx$
are normalized and converge to $y$ in norm. 

\endproof

\subsection{Properties of C-cones}

The operator $Y$ in Proposition \ref {prop:coneposop} may be non-unique.
Therefore the choice made in the proof of the proposition is made explicit.

\begin{definition}
\label{rel:def:Xy}
Given $y$ in ${\mathscr C}_x$ let $X_{y,x}$ denote the square of the Friedrichs extension of the
operator $Ax\mapsto Ay$, $A\in{\cal A}$.
\end{definition}

For further use note that the operator $X_{y,x}$ is affiliated with the commutant ${\cal A}'$ of $\cal A$
and that the vector $X_{y,x}^{t/2}x$ belongs to the cone ${\mathscr C}_x$ for any $t$ in $[0,1]$.

\begin{proposition}
\label{prop:prop:cyclic}
Take $y$ in ${\mathscr C}_x$. 
The following statements are equivalent.
\begin{description}
  \item [\quad (a)\quad] $y$ is cyclic for $\cal A$;
  \item [\quad (b)\quad] The operator $X^{1/2}_{y,x}$ is one-to-one;
  \item [\quad (c)\quad] $X^{t/2}_{y,x}x$ is separating for $\cal A$ for all $t\in [0,1]$;
  \item [\quad (d)\quad] $X^{t/2}_{y,x}x$ is cyclic
  for $\cal A$ for all $t\in [0,1]$.
\end{description}

\end{proposition}

\beginproof
\paragraph{(a) implies (b)}
Assume that $X^{1/2}_{y,x}u=0$. Then one has for all $A$ in $\cal A$
\be
0=(X_{y,x}^{1/2}u,Ax)=(u,X_{y,x}^{1/2}Ax) =(u,Ay).
\nnee
Since $y$ is cyclic by assumption this implies that $u=0$.
This shows that the linear operator $X^{1/2}_{y,x}$ is one-to-one.

\paragraph{(b) implies (c)}
Let $v=X_{y,x}^{t/2}x$. 
Assume $Av=0$ with $A$ in $\cal A$. Then one has
\be
0=Av=AX_{y,x}^{t/2}x=X_{y,x}^{t/2}Ax.
\nnee
Because $X_{y,x}$ is one-to-one it follows that $Ax=0$.
Because $x$ is separating for $\cal A$ one concludes that $A=0$. 
Hence, $v$ is separating for $\cal A$.

\paragraph{(c) implies (d)}
Let $v=X_{y,x}^{t/2}x$.
Without loss of generality assume $v$ is normalized. Then $\omega_v$ is a faithful state.
Hence there exists a vector $z$ in the natural positive cone $\mathscr P$ such that
$\omega_y=\omega_z$. See for instance Theorem 2.5.31 of \cite{BR79}.
The vector $z$ is separating for $\cal A$. It is well-known \cite{AH74} that
vectors in $\mathscr P$ are separating if and only if they are cyclic.
See for instance Proposition 2.5.30 of \cite{BR79}. Hence $z$ is cyclic.
This implies that the isometry $U$ in ${\cal A}'$ defined by 
\be
UAv=Az
\qquad A\in{\cal A},
\nnee
is actually a unitary operator.
Assume now that $(u,Av)=0$ for all $A$ in ${\cal A}$.
Then one has
\be
0&=&(u,Av)\cr
&=&(u,AUz)\cr
&=&(Uu,Az),
\qquad A\in{\cal A}.
\nnee
Because $z$ is cyclic for $\cal A$ it follows that $Uu=0$ and hence $u=0$.
This shows that $v$ is cyclic for $\cal A$.

\paragraph{(d) implies (a)}
This follows by taking $t=1$.

\endproof

\section{Exponential arcs in Hilbert space}
\label{sect:Hsp}

\subsection{Definition}

The notion of an {\em open exponential arc} was introduced in \cite {PC07,PG13}
and is used for instance in \cite{SST17}. Below is given a definition of
exponential arcs with values in Hilbert space. A small difference
in concept follows from referring to the arcs by their end points,
or by a midpoint and one end point,
rather than by two distinct points belonging to the interior.
In addition, uni-directional arcs are introduced
instead of the equivalence relations of \cite {PC07}.
The advantage is that they can be defined with any vector in ${\mathscr C}^1_x$
as the end point.

\begin{definition}
\label{exparchil:def:exparc}
Let  be given $y$ in ${\mathscr C}^1_x$ and let $X_{y,x}$ be the selfadjoint operator
affiliated with the commutant ${\cal A}'$ as defined in Definition \ref {rel:def:Xy}.
The {\em exponential arc} connecting $y$ to $x$ is the map $t\in [0,1]\mapsto \gamma_{y,x}(t)$
defined by
\be
\gamma_{y,x}(t)&=&e^{-\zeta_{y,x}(t)}\, X_{y,x}^{t/2}x, \qquad t\in [0,1],
\label{arc:def:arc}
\ee
where the normalization $\zeta_{y,x}(t)$ is given by
\be
\zeta_{y,x}(t)&=&\log ||X_{y,x}^{t/2}x||.
\nnee

\end{definition}

By the concavity of the function $u\mapsto u^t$ for $t$ in $[0,1]$ one has
\be
\zeta_{y,x}(t)=\log ||X_{y,x}^{t/2}x||\le t\log ||y||=0.
\nnee
Note that $\zeta_{y,x}(0)=\zeta_{y,x}(1)=0$ and $y_0=x$ and $y_1=y$.

For each $t$ in $[0,1]$ the vector $\gamma_{y,x}(t)$ belongs to ${\mathscr C}^1_x$.
It is properly normalized and one has for all positive $A$ in $\cal A$ that
\be
(\gamma_{y,x}(t),Ax)
&=&
e^{-\zeta_{y,x}(t)}\,\left(X_{y,x}^{t/2}x,Ax\right)\cr
&=&
e^{-\zeta_{y,x}(t)}\,\left(X_{y,x}^{t/4}x,AX_{y,x}^{t/4}x\right)\cr
&\ge&0.
\nnee

\begin{proposition}
The operator $X_{y,x}^{t/2}$ is the Friedrichs extension of the 
symmetric operator $Ax\mapsto AX_{y,x}^{t/2}x$ with $A$ in $\cal A$.
\end{proposition}

\def\tinyF{{\mbox{\tiny F}}}

\beginproof
Let $v=X_{y,x}^{t/2}x$ and let $Y$ denote the operator  
\be
Y:\,Ax\mapsto AX_{y,x}^{t/2}x, \qquad A\in{\cal A}.
\nnee
Introduce the sesquilinear form $s$ defined by
\be
s(Ax,Bx)&=&(Av,Bx),\qquad A,B\in{\cal A}.
\nnee
It is straightforward to verify that $s$ is well-defined, symmetric and positive-semi\-definite.
The domain of its closure $\overline s$ consists of all $u\in{\mathscr H}$ for which
there exist $u_n$ in ${\cal A}x$ converging to $u$ such that $s(u_n-u_m,u_n-u_m)$
tends to 0 as $n,m$ tend to $\infty$.
It is straightforward to verify that the domain of $X_{y,x}^{t/4}$, and a fortiori that of $X_{y,x}^{t/2}$,
is a subset of domain of $\overline s$.
Hence, both $Y^\tinyF$ and $X_{y,x}^{t/2}$ are selfadjoint extensions of the operator $Y$.
However, the Friedrichs extension $Y^\tinyF$ is the unique extension with domain included in $\dom \overline s$
--- see Theorem 2.11, p. 326 of \cite{KT66}. Hence, both operators coincide.

\endproof

\begin{corollary}
\label{exparchil:corol}
If $z=\gamma_{y,x}(t)$ with $0<t\le 1$ then $X^{1/2}_{z,x}=e^{-\zeta_{y,x}(t)}X^{t/2}_{y,x}$.
\end{corollary}

\begin{theorem}
\label{exparchil:prop:subarc}
Let $x$ be a normalized vector cyclic and separating for the von Neumann algebra $\cal A$.
For each $y$ in ${\mathscr C}^1_x$ and $0<t\le 1$ the map $s\mapsto \gamma_{y,x}(st)$
is the exponential arc connecting $\gamma_{y,x}(t)$ to $x$. The normalization functions satisfy
\be
\zeta_{y,x}(st)&=&\zeta_{z,x}(s)+s\zeta_{y,x}(t)
\quad\mbox{ with }\quad z=\gamma_{y,x}(t).
\label{exparc:thm:norm}
\ee
\end{theorem}

\beginproof
For convenience, let $z=\gamma_{y,x}(t)$.
We have to show that $\gamma_{y,x}(st)=v$ with $v=\gamma_{z,x}(s)$. 
The above corollary implies that
\be
X^{1/2}_{z,x}&=&e^{-\zeta_{y,x}(t)}X^{t/2}_{y,x}
\nnee
and
\be
X^{1/2}_{v,x}
&=&e^{-\zeta_{z,x}(s)}X^{s/2}_{z,x}\cr
&=&
e^{-\zeta_{z,x}(s)}e^{-s\zeta_{y,x}(t)}X^{st/2}_{y,x}
\nnee
and
\be
X^{1/2}_{v,x}
&=&
e^{-\zeta_{y,x}(st)}X^{st/2}_{y,x}.
\nnee
Comparison of the latter two gives (\ref {exparc:thm:norm}).
Application of the latter on the vector $x$ gives
\be
v&=&X^{1/2}_{v,x}x\cr
&=&e^{-\zeta_{y,x}(st)}X^{st/2}_{y,x}x\cr
&=&\gamma_{y,x}(st).
\nnee
\endproof

\subsection{The normalization function}
\label{sect:norm}

The following proposition mimics part of Proposition 6 of \cite{NJ19}.

\begin{proposition}
\label{norm:prop1}
Take  $y$ in ${\mathscr C}^1_x$. One has
\begin{description}
  \item [\quad (a) \quad] $0\le \zeta_{y,x}(st)\le s\zeta_{y,x}(t)$ for all $s,t$ in $[0,1]$;
  \item [\quad (b) \quad] The function $t\mapsto \zeta_{y,x}(t)$ is convex on the interval $[0,1]$;
  \item [\quad (c) \quad] The function $t\mapsto \zeta_{y,x}(t)$ is continuous at $t\downarrow 0$.
\end{description}
If $x$ belongs to the domain of $\log X_{y,x}$ then 
\begin{description}
 \item [\quad (d) \quad] The left derivative of $t\mapsto \zeta_{y,x}(t)$ at $t=0$ exists and one has
 \be
  \frac{\upd \,}{\upd t}\bigg|_{t=0}\zeta_{y,x}&=&\frac 12([\log X_{y,x}]x,x).
  \label{exparchil:prop:deriv}
  \ee
\end{description}

\end{proposition}

\beginproof
\paragraph{(a)}
This follows from (\ref {exparc:thm:norm}).

\paragraph{(b)}
Let
\be
X_{y,x}&=&\int_0^\infty\lambda \upd E_\lambda
\nnee
denote the spectral decomposition of the selfadjoint operator $X_{y,x}$.
Let $r=(1-s)t_1+st_2$ with $s$ in $[0,1]$ and $0\le t_1<t_2\le 1$.
Introduce the abbreviation
\be
\upd \mu(\lambda)&=&\lambda^{t_1}\upd(E_\lambda x,x).
\nnee
Without restriction of generality assume that the measure $\mu$ is non-trivial.
One has
\be
\int_0^\infty \lambda^r \upd (E_\lambda x,x)
&=&
\left[\int_0^\infty \upd \mu(\lambda)\right]\,
\frac{\int_0^\infty \lambda^{s(t_2-t_1)}\upd \mu(\lambda)}{\int_0^\infty \upd \mu(\lambda)}
\cr
&\le&
\left[\int_0^\infty \upd \mu(\lambda)\right]\,
\left[\frac{\int_0^\infty \lambda^{t_2-t_1}\upd \mu(\lambda)}{\int_0^\infty \upd \mu(\lambda)}
\right]^s\cr
&=&
\left[\int_0^\infty \upd \mu(\lambda)\right]^{1-s}\,
\left[\int_0^\infty \lambda^{t_2-t_1}\upd \mu(\lambda)\right]^s.
\nnee
Take the logarithm to obtain
\be
\zeta_{y,x}(r)
&=&\frac 12\log \int_0^\infty \lambda^r \upd (E_\lambda x,x)\cr
&\le&
\frac {1-s}2\log \left[\int_0^\infty \upd \mu(\lambda)\right]
+\frac s2\log \left[\int_0^\infty \lambda^{t_2-t_1}\upd \mu(\lambda)\right]\cr
&=&
(1-s)\zeta_{y,x}(t_1)+s\zeta_{y,x}(t_2).
\nnee
This shows the convexity of $t\mapsto \zeta_{y,x}(t)$.

\paragraph{(c)}
This follows from (a) in combination with (b).

\paragraph{(d)}
Starting point are the inequalities
\be
0\le e^t-1-t\le t^2,
\qquad -\frac 12\le t\le \frac 12.
\nnee
Let $0<\epsilon<1$.
One has for $0\le t\le 1/(-2\log\epsilon)$
\be
0\le
\int_{\epsilon}^{1/\epsilon}\left[\lambda^t-1-t\log\lambda\right]\upd(E_\lambda x,x)
&\le&
t^2\int_{\epsilon}^{1/\epsilon}
[\log\lambda]^2\upd(E_\lambda x,x).
\nnee
This can be written as
\be
\int_{\epsilon}^{1/\epsilon}[\log\lambda]\upd(E_\lambda x,x)
&\le&
\int_{\epsilon}^{1/\epsilon}\frac{\lambda^t-1}t\upd(E_\lambda x,x)\cr
&\le& 
\int_{\epsilon}^{1/\epsilon}[\log\lambda]\upd(E_\lambda x,x)
+t\int_{\epsilon}^{1/\epsilon}
[\log\lambda]^2\upd(E_\lambda x,x).
\nnee
Let $R(\epsilon)$ given by
\be
R(\epsilon)&=&
\left[\int_0^{\epsilon}+\int_{1/\epsilon}^{+\infty}\right]\,
\left[\log\lambda\right]\upd(E_\lambda x,x).
\nnee
One obtains
\be
-R(\epsilon)+([\log X_{y,x}]x,x)
&=&
\int_{\epsilon}^{1/\epsilon}[\log\lambda]\upd(E_\lambda x,x)\cr
&\le&
\int_{\epsilon}^{1/\epsilon}\frac{\lambda^t-1}t\upd(E_\lambda x,x)\cr
&\le&
-R(\epsilon)+([\log X_{y,x}]x,x)+t\,||[\log X_{y,x}]x||^2 
\nnee

Note that $R(\epsilon)$ tends to zero as $\epsilon\downarrow 0$ because by 
assumption $x$ belongs to the domain of $\log X_{y,x}$.
Hence, one concludes that
\be
\frac{\upd\,}{\upd t}\bigg|_{t=0}\,||X^{t/2}_{y,x}x||^2
&=&
([\log X_{y,x}]x,x).
\nnee
This leads to (\ref {exparchil:prop:deriv}).

\endproof

\subsection{Tangent vectors}

An important property of exponential arcs is differentiability.

\begin{proposition}
Given $y$ in ${\mathscr C}^1_x$ one has
\be
\dom X^{t/2}_{y,x}\log X_{y,x}&\supset&\dom X^{1/2}_{y,x},
\qquad 0<t<1.
\nnee
\end{proposition}

\beginproof
Let
\be
X_{y,x}&=&\int_0^\infty\lambda \upd E_\lambda
\nnee
denote the spectral decomposition of the selfadjoint operator $X_{y,x}$.
It follows from the spectral theorem that the domain of $X^{t/2}_{y,x}\log X_{y,x}$ is given by
\be
\dom X^{t/2}_{y,x}\log X_{y,x}
&=&
\{u:\,\int_0^{+\infty}\lambda^t(\log\lambda)^2\upd(E_\lambda u,u)<+\infty\}.
\nnee
Note that for $0\le\lambda<1$ the function $\lambda^t(\log\lambda)^2$ is bounded
while for large enough $\lambda$ one has $\lambda^t(\log\lambda)^2<\lambda$.
If $u$ belongs to $\dom X^{1/2}_{y,x}$ then one has
\be
\int_0^{+\infty}\lambda \upd(E_\lambda u,u)<+\infty
\nnee
so that $u$ belongs to $\dom X^{t/2}_{y,x}$.
\endproof

\begin{theorem}
\label{tangent:thm}
Let $x$ be a normalized vector cyclic and separating for the von Neumann algebra $\cal A$.
Take  $y$ in ${\mathscr C}^1_x$ and consider the exponential arc $t\mapsto \gamma_{y,x}(t)$
connecting  $y$ to $x$. The Fr\'echet derivative $\dot \gamma_{y,x}(t)$
of $t\mapsto \gamma_{y,x}(t)$ exists for $0< t<1$ and is given by
\be
\dot\gamma_{y,x}(t)
&=&
\left[\frac 12\log X_{y,x}- \zeta'_{y,x}(t)\right]\gamma_{y,x}(t),
\qquad 0<t<1,
\label{tan:thm:frech}
\ee
where $\zeta'_{y,x}(t)$ is the derivative of the normalization function $\zeta_{y,x}(t)$.
\end{theorem}

\beginproof
Let
\be
X_{y,x}&=&\int_0^\infty\lambda \upd E_\lambda
\nnee
denote the spectral decomposition of the selfadjoint operator $X_{y,x}$.
One has for $0<t<t+\epsilon<1$
\be
||X^{(t+\epsilon)/2}_{y,x}x-X^{t/2}_{y,x}x
-\frac \epsilon 2 [X^{t/2}_{y,x}\log X_{y,x}]x||^2
&=&
\int_0^{+\infty}f_t^2(\sqrt\lambda,\epsilon)\upd(E_\lambda x,x)
\nnee
with
\be
f_t(\mu,\epsilon)&=&\mu^{t+\epsilon}-\mu^{t}
-\epsilon \mu^{t}\log\mu.
\nnee
Note that $f_t(0,\epsilon)=f_t(1,\epsilon)=0$. Hence, $\mu\mapsto f_t^2(\mu,\epsilon)$
attains its maximum in the interval $[0,1]$ at some intermediate value $\overline\mu$ of $\mu$
that is the solution of the equation
\be
(t+\epsilon)\,\left[\mu^{\epsilon}-1\right]&=&\epsilon t\log\mu.
\nnee
In the limit of vanishing $\epsilon$ is
\be
\log\overline\mu&=&-\frac 2t-\frac{4\epsilon}{3t^3}+\mbox{O}^2(\epsilon).
\nnee
Use this to show that
\be
f_t(\overline\mu,\epsilon)
&=&\mbox{O}(\epsilon^2).
\nnee
One obtains
\be
\int_0^{1}f_t^2(\sqrt\lambda,\epsilon)\upd(E_\lambda x,x)
&\le&f_t^2(\overline\mu,\epsilon)\int_0^{1}\upd(E_\lambda x,x)\cr
&\le&f_t^2(\overline\mu,\epsilon)\cr
&=&\mbox{O}^2(\epsilon^4).
\label{tan:thm:temp1}
\ee

Next introduce the function $g_t(\nu,\epsilon)$ defined by
\be
f_t(\mu,\epsilon)=\sqrt\lambda\,g_t(\nu,\epsilon)
\quad\mbox{ with }\nu=1/\lambda.
\nnee
One finds $g_t(\nu,\epsilon)=f_{1-t}(\nu,-\epsilon)$.
Hence, the maximum with $\nu$ in the interval $[0,1]$ is reached 
for a value $\overline\nu$ satisfying
\be
g_t(\overline\mu,\epsilon)
&=&\mbox{O}(\epsilon^2).
\nnee
One has now
\be
\int_1^{+\infty}f^2(\sqrt\lambda,\epsilon)\upd(E_\lambda x,x)
&\le&
g^2(\overline\nu,\epsilon)\int_1^{+\infty}\lambda \upd(E_\lambda x,x)\cr
&\le&
g^2(\overline\nu,\epsilon)\,||X^{1/2}_{y,x}y||^2\cr
&=&
g^2(\overline\nu,\epsilon)\cr
&=&\mbox{O}(\epsilon^4).
\label{tan:thm:temp2}
\ee

Both estimates (\ref {tan:thm:temp1}) and (\ref {tan:thm:temp2})
together show that $(1/2) [X^{t/2}_{y,x}\log X_{y,x}]x$ is the Fr\'echet derivative
of the map $t\mapsto X^{t/2}_{y,x}x$. The existence of the Fr\'echet derivative
implies that also the normalization function $\zeta_{y,x}(t)$ is differentiable.
The result is (\ref {tan:thm:frech}).

\endproof

\subsection{Inverted arcs}

For convenience, introduce the following definition.

\begin{definition}
Given a vector $y$ in ${\mathscr C}^1_x$ 
the set $\mathring\gamma_{y,x}$ is the image of $[0,1]$ by the exponential arc 
$t\in [0,1]\mapsto \gamma_{y,x}(t)$ minus the endpoints $x$ and $y$. 
\be
\mathring\gamma_{y,x}&=&\{\gamma_{y,x}(t):\, t\in[0,1]\}\setminus\{x,y\}\subset{\mathscr H}.
\nnee
\end{definition}

\begin{proposition}
\label{inv:prop:allcs}
If $y$ in ${\mathscr C}^1_x$ is cyclic for $\cal A$ then all vectors
$\gamma_{y,x}(t)$ are cyclic and separating for $\cal A$.
\end{proposition}

\beginproof
It follows from Proposition \ref{prop:prop:cyclic} that $\gamma_{y,x}(t)$ is a cyclic vector. 
It is separating by Proposition \ref {prop:prop:cyclic}.
\endproof

If $y$ in ${\mathscr C}^1_x$ is cyclic for $\cal A$ then
the role of the vectors $x$ and $y$ may be interchanged and the arc can be inverted.

\begin{proposition}
\label{inv:prop:equiv}
If $y$ in ${\mathscr C}^1_x$ is cyclic for $\cal A$ then $x$ belongs to ${\mathscr C}^1_y$
and $t\mapsto \gamma_{y,x}(1-t)$ is the exponential arc connecting $x$ to $y$.
The normalization functions satisfy $\zeta_{x,y}(t)=\zeta_{y,x}(1-t)$.
The sets $\mathring\gamma_{y,x}$ and $\mathring\gamma_{x,y}$ coincide.
\end{proposition}

\beginproof
If $y$ is cyclic then the operator $X_{y,x}$ has a densely defined inverse.
The vector $x$ belongs to ${\mathscr C}^1_y$ because it can be written as 
$x=X^{-1/2}_{y,x}y$ with $X^{-1/2}_{y,x}$ a selfadjoint operator affiliated with
the commutant ${\cal A}'$. The inverse of the Friedrichs extension of the map
$Ax\mapsto Ay$ is the Friedrichs extension of the map $Ay\mapsto Ax$. Hence one has 
$X^{1/2}_{x,y}=X^{-1/2}_{y,x}$.

The remainder of the proof is then straightforward and follows from
\be
X^{(1-t)/2}_{y,x}x=X^{-t/2}_{y,x}y=X^{t/2}_{x,y}y.
\nnee

\endproof

The following result extends that of Theorem \ref {exparchil:prop:subarc}
in the case that all vectors of the exponential arc are cyclic for the von
Neumann algebra $\cal A$. To do so it uses Proposition \ref {inv:prop:equiv}.

\begin{theorem}
\label{exparcstate:prop:subarc}
Let $x$ be a normalized vector cyclic and separating for the von Neumann algebra $\cal A$.
Assume $y$ in ${\mathscr C}^1_x$ is cyclic for $\cal A$. Then it holds that
for any $s$ and $t$ in $[0,1]$ the map $r\mapsto \gamma_{y,x}((1-r)s+rt)$
is the exponential arc connecting  $\gamma_{y,x}(t)$ to $\gamma_{y,x}(s)$.
\end{theorem}

\beginproof
Proposition \ref {inv:prop:equiv} allows us without loss of generality 
to make the assumption that $0<s<t\le 1$.

For convenience, let $z=\gamma_{y,x}(s)$ and $u=\gamma_{y,x}(t)$.
We have to show that $\gamma_{u,z}(r)=\gamma_{y,x}((1-r)s+rt)$.

Let $q=r(t-s)/(1-s)$.
It follows from Theorem \ref {exparchil:prop:subarc} that
$\gamma_{u,z}(r)=\gamma_{y,z}(q)$.
Proposition \ref {inv:prop:equiv} yields 
$\gamma_{y,z}(q)=\gamma_{z,y}(1-q)$.
Finally, apply Theorem \ref {exparchil:prop:subarc} again to obtain 
$\gamma_{z,y}(1-q)=\gamma_{x,y}(rs+(1-r)t)$ and 
Proposition \ref {inv:prop:equiv} to obtain 
$\gamma_{x,y}(rs+(1-r)t)=\gamma_{y,x}((1-r)s+rt)$. 
\endproof

\subsection{Midpoints of exponential arcs}

Theorem \ref {tangent:thm} asserts that any exponential arc $t\mapsto \gamma_{y,x}(t)$
is Fr\'echet differentiable on the open interval $(0,1)$. The end points $t=0$ and 
$t=1$ are not included. This suggests to consider exponential arcs that have $x$
as a midpoint.

\begin{definition}
Let $\mathring\gamma_x$ denote the union of all sets $\mathring\gamma_{z,y}$
where $y$ and $z$ are vectors in ${\mathscr C}^1_x$
that are cyclic and separating for $\cal A$ with the property that $x=\gamma_{z,y}(1/2)$.
\end{definition}

It follows from Theorem \ref {exparchil:prop:subarc} that any vector in $\mathring\gamma_{z,y}$
with $y$ and $z$ in ${\mathscr C}^1_x$ and with midpoint $x$ also belongs to ${\mathscr C}^1_x$.
Hence one has the inclusion $\mathring\gamma_x\subset {\mathscr C}^1_x$.
Vectors $y$ in ${\mathscr C}^1_x$ are excluded from $\mathring\gamma_x$
if they are not cyclic for $\cal A$ or the exponential arc $t\mapsto\gamma_{y,x}(t)$
cannot be extended beyond the initial vector $x$ or the final vector $y$.

\begin{proposition}
\label{mid:prop:extend}
For any $u$ in $\mathring\gamma_x$, $u\not=x$,
the exponential arc $t\in[0,1]\mapsto \gamma_{u,x}(t)$
can be extended to a map $t\in [-1,1]\mapsto \gamma_{u,x}(t)$ with the following properties.
\begin{description}
 \item [\quad (a) \quad] There exist $s$ in $(1/2,1]$ and normalized vectors $y$ and $z$
        that are cyclic and separating for $\cal A$ such that
        \be
            \gamma_{u,x}(t)&=&\gamma_{z,y}((s-1/2)t+1/2),
            \qquad -1\le t\le 1;
            \label{mid:extexparc}
        \ee
 \item [\quad (b) \quad] The Fr\'echet derivative $\dot \gamma_{u,x}(t)$ of 
        $\gamma_{u,x}(t)$ exists for all $t$ in the open interval $(-1,1)$;
        it is given by
        \be
            \dot \gamma_{u,x}(t)
            &=&
            (s-1/2)\dot\gamma_{z,y}(r)
            \label{mid:prop:frech}
        \ee
        with $r=(s-1/2)t+1/2$;
 \item [\quad (c) \quad] The Fr\'echet derivative $\dot \gamma_{u,x}(t)$ satisfies the equation
 \be
\dot\gamma_{u,x}(t)
&=&
\left[\frac 12\log X_{u,x}- \zeta'_{u,x}(t)\right]\gamma_{u,x}(t),
\qquad -1<t<1,
\label{mid:thm:frech}
\ee
with $\zeta_{u,x}(t)=\log||X^{t/2}_{u,x}x||$.
\end{description}

\end{proposition}

\beginproof

\paragraph{(a)}
From the definition of $\mathring\gamma_x$ it follows that there exist 
normalized vectors $y$ and $z$
that are cyclic and separating for $\cal A$ with the properties that $x=\gamma_{z,y}(1/2)$
and $u=\gamma_{z,y}(s)$ for some $s$ in $(1/2,1]$.
It then follows from Theorem \ref {exparchil:prop:subarc}
that $\gamma_{u,x}(t)=\gamma_{z,y}(r)$ holds for $t$ in $[0,1]$.
Hence, the definition (\ref {mid:extexparc}) is an extension of the exponential arc $\gamma_{u,x}$
from $[0,1]$ to $[-1,1]$.

\paragraph{(b)}
The existence of the Fr\'echet derivative is proved in Theorem \ref {tangent:thm}.
It satisfies
\be
\dot\gamma_{z,y}(r)
&=&
\left[
\frac 12 \log X_{z,y}- \zeta'_{z,y}(r)
\right]\gamma_{z,y}(r).
\label{mid:prop:temp}
\ee
This gives (\ref {mid:prop:frech}).

\paragraph{(c)}

Apply Corollary \ref {exparchil:corol} with $u=\gamma_{z,y}(s)$ and $x=\gamma_{z,y}(1/2)$
to obtain
\be
X_{z,y}^{s/2}&=&e^{\zeta_{z,y}(s)}X_{u,y}^{1/2},\cr
X_{z,y}^{1/4}&=&e^{\zeta_{z,y}(1/2)}X^{1/2}_{x,y}.
\ee
Together these expressions yield
\be
\frac 12(s-1/2)\log X_{z,y}
&=&
\zeta_{z,y}(s)+\frac 12\log X_{u,y}-\zeta_{z,y}(1/2)-\frac 12\log X_{x,y}.
\nnee
Hence one has
\be
\dot\gamma_{u,x}(t)
&=&
\left(s-\frac 12\right)\dot\gamma_{z,y}(r)\cr
&=&
\left(s-\frac 12\right)
\left[\frac 12[\log X_{z,y}]-\zeta'_{z,y}(r)\right]\gamma_{u,x}(t)\cr
&=&
\left[\frac 12\log X_{u,y}-\frac 12\log X_{x,y}+\zeta_{z,y}(s)-\zeta_{z,y}(1/2)-\zeta'_{z,y}(r)\right]
\gamma_{u,x}(t).
\nnee
Comparison with (\ref {tan:thm:frech}) gives
\be
[\log X_{u,x}]\gamma_{u,x}(t)
&=&
\left[\log X_{u,y}-\log X_{x,y}+f(t)\right]\gamma_{u,x}(t),
\qquad 0<t<1,
\nnee
for some real function $f(t)$.

\endproof

\section{Connecting vector states}
\label{sect:ssp}

\subsection{Definitions and properties}
\label{sect:connect}

The exponential arcs of vectors in Hilbert space define exponential arcs of vector states
on the von Neumann algebra $\cal A$.

\begin{definition}
Let $\Vo_x$ denote the subset of $\Vo$ consisting of the vector states $\omega_y$
where $y$ belongs to ${\mathscr C}^1_x$.
\end{definition}

\begin{proposition}
\label{defcon:prop:equiv}
Let $x$ and $y$ be normalized vectors in $\mathscr H$ that are cyclic and separating for $\cal A$.
If $\omega_x=\omega_y$ then there exists a unitary operator $U$ in the commutant algebra ${\cal A}'$
which maps ${\mathscr C}^1_x$ onto ${\mathscr C}^1_y$ and 
exponential arcs onto exponential arcs in such a way that
\be
\gamma_{Uz,y}(t)&=&U\gamma_{z,x}(t),
\qquad z\in {\mathscr C}^1_x
\mbox{ and }t\in [0,1].
\label{defcon:prop:map}
\ee
In particular, one has $\Vo_x=\Vo_y$. 
\end{proposition}

\beginproof
A unitary operator $U$ is defined by $UAx=Ay$ for all $A$ in $\cal A$.
It belongs to the commutant algebra ${\cal A}'$.
Take now $z$ in ${\mathscr C}^1_x$. Then one has
\be
(Uz,Ay)=(Uz,UAx)=(z,Ax)\ge 0,
\qquad A\in{\cal A}, A\ge 0.
\nnee
This shows that $Uz$ belongs to ${\mathscr C}^1_y$.
Now one has
\be
\omega_{Uz}(A)=(AUz,Uz)=(UAz,Uz)=(Az,z)=\omega_z(A),
\qquad A\in {\cal A}.
\nnee
This shows that $\omega_z$ belongs to $\Vo_y$ and hence that $\Vo_x\subset\Vo_y$.
Equality follows by interchanging $x$ and $y$.

The proof of (\ref {defcon:prop:map}) is straightforward.

\endproof

The above proposition justifies the following definition.

\begin{definition}
Given $y$ in ${\mathscr C}^1_x$ a map $t\in [0,1]\mapsto \omega_t$ is an exponential arc
connecting the vector state $\omega_y$ to the vector state $\omega_x$ if there exists
an exponential arc $t\in [0,1]\mapsto \gamma_{y,x}(t)$ connecting $y$ to $x$ such that
$\omega_t=\omega_{y_t}$ with $y_t={\gamma_{y,x}(t)}$ for all $t$ in $[0,1]$.
\end{definition}

Because the map $y\in{\mathscr C}^1_x\mapsto \omega_y$ is one-to-one there is only one such
exponential arc for each $y$ in ${\mathscr C}^1_x$. 
In addition, it follows immediately from Proposition \ref {exparchil:prop:subarc}
that for each $t$ satisfying $0<t\le 1$ the map $s\mapsto \omega_{y_{st}}$ is the unique
exponential arc connecting $\omega_{y_t}$ to $\omega_x$.
Proposition \ref {defcon:prop:equiv} shows that the exponential arc does not depend
on the choice of vector $x$ representing the state $\omega_x$.

In summary, one has the following.

\begin{proposition}
Any state $\omega$ of $\Vo_x$ is connected to the state $\omega_x$ by a unique exponential arc
$t\mapsto\omega_t$
such that $\omega_0=\omega_x$ and $\omega_1=\omega_y$.
\end{proposition}

If $t\mapsto\omega_t$ is an exponential arc connecting $\omega_y$ to $\omega_x$ 
and the normalization function $\zeta_{y,x}(t)$ is strictly convex
then the set of vector states $\{\omega_{y_t}:\, 0<t<1\}$
is an exponential model.
The primary chart is the map $\omega_{y_t}\mapsto t$. The tangent vectors $\dot\omega_{y_t}$ belong to the
predual ${\cal A}_*$. This is a subspace of the dual ${\cal A}^*$ of $\cal A$
and is a Banach space for the norm of ${\cal A}^*$. See for instance
Section 2.4.3 of \cite{BR79}. These tangent vectors are given by
\be
\dot\omega_{y_t}(A)
&=&
(A\dot\gamma_{y,x}(t),\gamma_{y,x}(t))+(A\gamma_{y,x}(t),\dot \gamma_{y,x}(t))\cr
&=&
e^{-\zeta_{y,x}(t)}(AX^{t/2}_{y,x}x,[\log X_{y,x}]X^{t/2}_{y,x}x)
-2\zeta'_{y,x}(t)\omega_{y_t}(A),
\qquad A\in{\cal A}.
\nnee
The Legendre transformation of $\zeta_{y,x}(t)$ is given by
\be
\zeta^*_{y,x}(s)
&\equiv&
\sup\{st-\zeta_{y,x}(t):\,0\le t\le 1\}.
\nnee
If $\zeta_{y,x}(t)$ is strictly convex then the supremum is reached at a unique value of $s$.
The latter is denoted $t^*$ and is given by $t^*=\zeta'_{y,x}(t)$. The following identity holds
\be
\zeta_{y,x}(t) +\zeta^*_{y,x}(t^*)&=&tt^*.
\nnee
The map $\omega_t \mapsto t^*=\zeta'_{y,x}(t)$ is the dual chart of the model.

\subsection{Majorizing states}

\begin{definition}
Fix a normalized vector $x$ in $\mathscr H$ that is cyclic and separating for $\cal A$.
Let $\breve{\Vo}_x$ denote the subset of states of ${\Vo}$ that are
majorized by a multiple of the state $\omega_x$.
\end{definition}

Additional properties can be proved for vector states belonging to  $\breve{\Vo}_x$.
See Section 9 of \cite{AH74}. The following is an adaptation of some of these results
to the present context.

\begin{proposition}
\label{major:prop}
Let $\omega\in\Vo$.
Are equivalent:
\begin{description}
  \item [\quad (a) \quad] $\omega$ belongs to $\breve{\Vo}_x$;
  \item [\quad (b) \quad] There exists $y$ in ${\mathscr C}^1_x$ such that $\omega=\omega_y$ 
  and the operator $X_{y,x}$ is bounded.
\end{description}

\end{proposition}

\beginproof

\paragraph{(a) implies (b)}

By assumption there exists $z$ such that $\omega=\omega_z$
and there exists a real constant $\lambda$ such that
\be
\omega_z(A)\le\lambda\omega_x(A)
\quad\mbox{ for all positive }A\in {\cal A}.
\nnee
The linear operator $Z$ defined by
\be
Z:\, Ax\mapsto Az,
\qquad A\in{\cal A},
\nnee
is well-defined because $x$ is cyclic and separating for $\cal A$.
It extends by continuity to a bounded operator in ${\cal A}'$ because
\be
||ZAx||^2=||Az||^2=\omega_z(A^*A)\le\lambda\omega_x(A^*A)=\lambda||Ax||^2,
\qquad A\in{\cal A}.
\nnee
Let $Z=J|Z|$ be the polar decomposition of $Z$ and let $Y=|Z|$ and $y=Yx$.
Then one has $\omega_y=\omega_z$ and 
\be
(y,A^*Ax)= ||A|Y|^{1/2}x||^2\ge 0,
\qquad A\in{\cal A}.
\nnee
This shows that $y$ belongs to ${\mathscr C}^1_x$ and hence that $\omega$ belongs to
${\Vo}_x$. Finally note that $X_{y,x}=Y^2=Z^*Z$.

\paragraph{(b) implies (a)}
One has for all $A$ in $\cal A$
\be
\omega_y(A^*A)&=&(AX^{1/2}_{y,x}x,AX^{1/2}_{y,x}x)=(X^{1/2}_{y,x}Ax,X^{1/2}_{y,x}Ax)\cr
&\le&
||X_{y,x}||\,(Ax,Ax)=||X_{y,x}||\,\omega_x(A^*A).
\nnee
This shows that $\omega_y$ is majorized by $||X_{y,x}||\,\omega_x$.

\endproof

\begin{corollary}
$\breve\Vo_x\subset \Vo_x$.
\end{corollary}

If the operator $X_{y,x}$ is bounded then one does not need to take care about domain problems.
From Definition \ref {exparchil:def:exparc} one obtains for $t\in [0,1]$
\be
\gamma_{y,x}(t)&=&e^{-\zeta_{y,x}(t)}\, \exp\left(-\frac t2H_{y,x}\right)x,
\nnee
with
\be
H_{y,x}&=&-\log X_{y,x}.
\nnee
The normalization $\zeta_{y,x}(t)$ can be written as
\be
\zeta_{y,x}(t)&=&\frac 12\log\left(\exp\left(-tH_{y,x}\right)x,x\right).
\nnee

The following result follows immediately from Proposition \ref {cones:prop:dense}.

\begin{proposition}
$\breve{\Vo}_x$ is norm-dense in ${\Vo}_x$.
\end{proposition}

\begin{proposition}
One has the following.
\begin{description}
  \item [\quad (a) \quad] The set of states $\breve{\Vo}_x$ is convex;
  \item [\quad (b) \quad] The norm closure $\overline {\Vo}_x$
  of the set of states ${\Vo}_x$ is convex.
\end{description}

\end{proposition}

\beginproof
\paragraph{(a)}
Take $y$ and $z$ in ${\mathscr C}^1_x$ and $\lambda$ in $[0,1]$.
If $\omega_y$ and $\omega_z$ belong to $\breve{\Vo}_x$
then they are dominated by multiples of $\omega_x$ and so is
$\lambda \omega_y+(1-\lambda)\omega_z$. It then follows from Proposition \ref {major:prop}
that the convex combination belongs to $\breve{\Vo}_x$.

\paragraph{(b)}
Take $\omega^{(1)}$ and $\omega^{(2)}$ in $\overline {\Vo}_x$ and let $\lambda\in [0,1]$.
There exist states $\omega^{(1)}_n$ and $\omega^{(2)}_n$ in  $\breve{\Vo}_x$
converging to $\omega^{(1)}$, respectively $\omega^{(2)}$.
The states $\lambda \omega^{(1)}_n+(1-\lambda)\omega^{(2)}_n$ belong to $\breve{\Vo}_x$
by item (a) of the proposition
and converge to $\lambda \omega^{(1)}+(1-\lambda)\omega^{(2)}$.
Hence, the latter state belongs to the closure of $\breve{\Vo}_x$
which coincides with the closure of ${\Vo}_x$.

\endproof

\begin{proposition}
If $y$ belongs to ${\mathscr C}^1_x$ and is cyclic for $\cal A$
   then one has  $\breve{\Vo}_x=\breve{\Vo}_y$ and $\overline {\Vo}_x=\overline {\Vo}_y$.
\end{proposition}

\beginproof
If $y$ belongs to ${\mathscr C}^1_x$ and 
is cyclic for $\cal A$ then Proposition \ref {inv:prop:equiv} asserts that
$x$ belongs to ${\mathscr C}^1_y$ as well.
Hence there exist positive constants $\lambda_{y,x}$ and $\lambda_{x,y}$ such that
\be
\omega_y\le\lambda_{y,x}\omega_x
\quad\mbox{ and }\quad
\omega_x\le\lambda_{x,y}\omega_y
\nnee
Take now $z$ in $\breve{\Vo}_x$. Then there exists $\lambda_{z,x}$ such that
\be
\omega_z\le\lambda_{z,x}\omega_x\le \lambda_{z,x}\lambda_{x,y}\omega_y.
\nnee
Hence, $\omega_z$ belongs to $\breve{\Vo}_y$. 
This shows that $\breve{\Vo}_x\subset\breve{\Vo}_y$.
By symmetry one concludes that $\breve{\Vo}_x=\breve{\Vo}_y$.
Because they are dense in $\overline {\Vo}_x$ respectively $\overline{\Vo}_y$
one concludes that also $\overline {\Vo}_x=\overline {\Vo}_y$.

\endproof

\subsection{A neighborhood of faithful states}

\begin{definition}
The set of vector states $\omega_y$ with $y$ in $\mathring\gamma_x$ is denoted
$\mathring\Vo_x$.
\end{definition}

\begin{proposition}
Take $y$ in ${\mathscr C}^1_x$. If $X_{y,x}$ is bounded with bounded inverse then
all vector states $\omega_u$ with $u=\gamma_{y,x}(t)$ and $-1<t<1$ belong to $\mathring\Vo_x$.
\end{proposition}

\beginproof
Let $z=X_{y,x}^{-1/2}$. Then $x$ is the midpoint of the exponential arcs $\gamma_{z,y}(t)$
and $\gamma_{y,z}(t)$. By Proposition \ref {inv:prop:allcs} all vectors of these arcs are 
cyclic and separating for $\cal A$. Hence, one has
$\mathring\gamma_{y,x}\subset \mathring\gamma_x$. From the above definition it then follows
that for all $z$ in $\mathring\gamma_{y,x}$ the vector state $\omega_z$ belongs to $\mathring\Vo_x$.

\endproof

The bounded operators with bounded inverse and belonging to the commutant ${\mathscr A}'$
form a real Banach-Lie group in the topology induced by the norm operator topology.
This observation is the starting point of \cite {CIJM19,CIJM19a}.

\begin{theorem}
\label{neigh:thm}
Let $x$ be a normalized vector cyclic and separating for the von Neumann algebra $\cal A$.
Choose a vector state $\omega$ in $\mathring\Vo_x$. 
Let $t\in [0,1]\mapsto\omega_t$ be the exponential arc connecting $\omega$ to $\omega_x$.
One has the following.
\begin{description}
 \item [\quad (a) \quad] The exponential arc $t\in[0,1] \mapsto\omega_t$ 
 can be extended to the interval $[-1,1]$ in such a way that $t\mapsto \omega_{2t-1}$
 is the exponential arc connecting $\omega_1$ to $\omega_{-1}$;
 \item [\quad (b) \quad] The Fr\'echet derivative $\dot\omega_t$ of $t\mapsto\omega_t$
 exists for $t$ in the open interval $(-1,1)$ and is given by
 \be
   \dot\omega_t(A)&=&(Ay_t,[\log X_{y,x}]y_t)-(Ay_t,y_t)(y_t,[\log X_{y,x}]y_t),
 \nnee 
 where $y_t$ is the unique vector in ${\mathscr C}_x$ such that $\omega_t=\omega_{y_t}$.
 The operator $X_{y_1,x}$ is the positive operator affiliated with
 the commutant ${\cal A}'$ as given by Definition \ref {rel:def:Xy}.
\end{description}
\end{theorem}

\beginproof

\paragraph{(a)}
From the definition of the exponential arc $t\in[0,1] \mapsto\omega_t$
it follows that there exist $y$ in ${\mathscr C}^1_x$ such that $\omega_t=\omega_{y_t}$
with $y_t=\gamma_{y,x}(t)$.
Because $\omega=\omega_1=\omega_{y}$ belongs to $\mathring\Vo_x$ there exists $u$ in $\mathring\gamma_x$
such that $\omega=\omega_u$. From the definition of $\mathring\gamma_x$
it follows that there exist $y,z$ in ${\mathscr C}^1_x$ such that $u\in \mathring\gamma_{y,z}$,
$x$ is the midpoint of $\gamma_{y,z}$ and $y$ and $z$ are cyclic and separating for $\cal A$.
By Theorem \ref {exparcstate:prop:subarc} $x$ is also the midpoint of $\gamma_{u,v}$ for some
$v$ belonging to the exponential arc $\gamma_{y,z}$.
By Proposition \ref {mid:prop:extend} the exponential arc connecting $u$ to $x$ can now be extended
in such a way that $\gamma_{u,x}(-1)=v$. Now $\omega_t$ can be extended to all of $[-1,1]$ by
putting $\omega_t=\omega_{y_t}$ with $y_t=\gamma_{u,x}(t)$.
By construction it is the exponential arc connecting $u$ to $v$.

\paragraph{(b)}
One calculates with the help of 
\be
\dot\omega_t(A)
&=&
\frac{\upd\,}{\upd t}(Ay_t,y_t)\cr
&=&
(A\dot y_t,y_t)+(Ay_t,\dot y_t)
\nnee
with
\be
\dot y_t&=&\left[\frac 12\log X_{y,x}-\zeta'_{y,x}(t)\right]y_t.
\nnee
It satisfies
\be
& &
|\omega_{t+\epsilon}(A)-\omega_t(A)-\epsilon\dot\omega_t(A)|\cr
&=&
|(Ay_{t+\epsilon},y_{t+\epsilon})-(Ay_t,y_t)-\epsilon(Ay_t,\dot y_t)-\epsilon(A\dot y_t,y_t)|\cr
&\le&
|(Ay_{t+\epsilon},y_{t+\epsilon})-(Ay_{t+\epsilon},y_t)-\epsilon(Ay_t,\dot y_t)|\cr
& &
+|(Ay_{t+\epsilon},y_t)-(Ay_t,y_t)-\epsilon(A\dot y_t,y_t)|\cr
&\le&
||A||\,||y_{t+\epsilon}-y_t-\epsilon \dot y_t)||\cr
& &
+||A||\,||y_{t+\epsilon},y_t-Ay_t-\epsilon \dot y_t||\cr
&=&\mbox{ o}(\epsilon).
\nnee
This shows that $\dot\omega_t$ is the Fr\'echet derivative of $t\mapsto\omega_t$.

\endproof

\section{Special cases}
\label{sect:special}

\subsection{The finite-dimensional commutative case}

Let $\cal A$ be the algebra of diagonal $n$-by-$n$ matrices with complex entries.
They act on the Hilbert space $\Co^n$. Let $x$ be a normalized element of $\Co^n$
with non-vanishing components. The corresponding probabilities are denoted $p_i$
and are given by $p_i=|x_i|^2$. The state $\omega_x$ on $\cal A$ is defined by
\be
\omega_x(A)=(Ax,x)=\sum_{i=1}^np_i A_{i,i}. 
\nnee

The cone ${\mathscr C}_x$ consists of element $y$ of $\Co^n$ satisfying $y_i\overline {x_i}\ge 0$
for $i=1,2,\cdots n$. 
Hence, for each $y$ in ${\mathscr C}_x$  there exists a non-negative number $\tilde y_i$ such that
\be
y_i=\tilde y_i\,x_i.
\nnee
In particular, the operator $X_{y,x}$ is the diagonal
matrix with entries $(X_{y,x})_{i,i}=(\tilde y_i)^2$.
One verifies that for any $A$ in $\cal A$ one has
\be
(A X^{1/2}_{y,x}x,X^{1/2}_{y,x}x)
&=&
\sum_iA_{i,i}\,(X_{y,x})_{i,i}|x_i|^2\cr
&=&
\sum_iA_{i,i}(\tilde y_i)^2|x_i|^2\cr
&=&
\sum_iA_{i,i}|y_i|^2\cr
&=&
\omega_y(A).
\nnee

Note that $\breve{\mathscr C}_x={\mathscr C}_x$ because all operators on $\Co^n$ are bounded.
The algebra $\cal A$ coincides with its commutant ${\cal A}'$.
Hence, cone and dual cone coincide.  In this commutative example all cones  ${\mathscr C}_x$
are self-dual. The natural positive cone consists of all vectors $x$ with non-negative
components $x_i\ge 0$.

\paragraph{In Hilbert space}

The exponential arc $\gamma_{y,x}(t)$ is given by (\ref {arc:def:arc}).
This gives
\be
(\gamma_{y,x}(t))_i
&=&
e^{-\zeta_{y,x}(t)}\, (X_{y,x}^{t/2})_{i,i}x_i\cr
&=&
e^{-\zeta_{y,x}(t)}\, (\tilde y_i)^{t}x_i.
\nnee
Note that $y_i=0$ implies that $y_i^{t}=0$ for any $t>0$.
However, if $t=0$ then by convention $(\tilde y_i)^t=1$ whatever the value
is of $\tilde y_i$. Hence, discontinuities can occur at $t=0$.

The normalization function is given by
\be
\zeta_{y,x}(t)
&=&
\log ||X_{y,x}^{t/2}x||\cr
&=&
\frac 12\log\sum_i (\tilde y_i)^{2t}|x_i|^2.
\nnee
If one or more of the $y_i$ vanish then the function $\zeta_{y,x}(t)$ is discontinuous at $t=0$.
If $t>0$ then the function is differentiable and one has 
\be
\zeta'_{y,x}(t)
&=&
e^{-2\zeta_{y,x}(t)}\sum_i\strut'(\log\tilde y_i)(\tilde y_i)^{2t}|x_i|^2,
\nnee
where the sum $\Sigma'$ is restricted to those terms for which $\tilde y_i\not=0$.
The second derivative then becomes
\be
\zeta''_{y,x}(t)
&=&
2e^{-2\zeta_{y,x}(t)}\sum_i\strut'(\log\tilde y_i)^2(\tilde y_i)^{2t}|x_i|^2
-2\left[\zeta'_{y,x}(t)\right]^2\cr
&=&
2e^{-2\zeta_{y,x}(t)}\sum_i\strut'\left[\log\tilde y_i-\zeta'_{y,x}(t)\right]^2\,(\tilde y_i)^{2t}|x_i|^2.
\nnee

The vector $x$ is a midpoint of an arc ending at a vector $y$ of ${\mathscr C}^1_x$
if and only if $y_i\not=0$ for all $i$. If this is the case then by Proposition \ref {mid:prop:extend}
the path $t\mapsto\gamma_{y,x}(t)$
can be extended to the interval $[-1,1]$. It is actually clear from the explicit expression
given above that it can be extended to all of the real line.  
The set $\mathring\gamma_x$ consists of all $y$ in ${\mathscr C}^1_x$ such that
$y_i\not=0$ for all $i$. In particular, it does not depend on the choice of $x$.

\paragraph{In state space}

Starting from the assumption that $x_i\not=0$ for all $i$ one can
easily show that any state on the algebra $\cal A$ of diagonal matrices is majorized by a
multiple of the state $\omega_x$.
Hence, by Proposition \ref {major:prop} all states on $\cal A$ belong to the neighborhood
$\breve\Vo_x$ and {\em a fortiori} to $\Vo_x$.

Take $y$ in ${\mathscr C}^1_x$.
Let $p_i=|x_i|^2$ and $q_i=|y_i|^2$. The exponential arc connecting $\omega_y$ to $\omega_x$
satisfies the equation
\be
\dot\omega_t(A)&=&\omega_t(AH)-\omega_t(A)\omega_t(H)
\nnee
with the generator $H$ given by
\be
H_{i,i}=[\log X_{y,x}]_{i,i}=\log\frac{q_i}{p_i}.
\nnee
In particular the tangent vector at $t=0$ is given by
\be
\dot\omega_t(A)\bigg|_{t=0}&=&
\omega_x(AH)+D(\omega_x||\omega_y)\omega_x(A)
\nnee
with the divergence $D(\omega_x||\omega_y)$ given by
\be
D(\omega_x||\omega_y)&=&\sum_ip_i\log\frac{p_i}{q_i}.
\nnee

\subsection{The density matrix}

What follows here is partly based on \cite{NJ18,NJ19b}.

Consider the von Neumann algebra of $n$-by-$n$ matrices.
In its standard form
it acts on the Hilbert space ${\mathscr H}=\Co^n\otimes\Co^n$ (up to a $*$-isomorphism). 
Any vector $x$ of ${\mathscr H}$ can be written as a linear combination of vectors of the form $x^{(1)}\otimes x^{(2)}$.
Any operator $A$ of $\cal A$ is of the form $A=A^{(1)}\otimes\Io$, where
$A^{(1)}$ is an $n$-by-$n$ matrix and $\Io$ is the identity matrix.
Any element $A'$ of the commutant ${\cal A}'$ is of the form $A'=\Io\otimes A^{(2)}$.

For each state $\omega$ on $\cal A$ there exists
a unique density matrix $\rho_\omega$ such that 
\be
\omega(A)=\Tr\rho_\omega A^{(1)}.
\nnee
Let $(e_i)_i$ be a basis of eigenvectors of $\rho_\omega$
with corresponding eigenvalues $p_i$. Then a vector $x$ is given by
\be
x&=&\sum_i\sqrt{p_i}e_i\otimes e_i.
\label{dens:xvec}
\ee
It satisfies
\be
\omega_x(A)
&=&
\sum_{i,j}\sqrt{p_ip_j}\left(A^{(1)}e_i\otimes e_i,e_j\otimes e_j\right)\cr
&=&
\sum_ip_i\left(A^{(1)}e_i,e_i\right)\cr
&=&
\Tr\rho_\omega A^{(1)}.
\nnee

Assume that $p_i>0$ for all $i$. Then the vector $x$, as given by (\ref {dens:xvec}), is 
cyclic and separating for $\cal A$.
A straightforward calculation shows that the corresponding modular operator $\Delta_x$
equals $\rho_\omega\otimes\rho_\omega^{-1}$.

Let $\omega_y$ be an arbitrary vector state with $y$ in ${\mathscr C}^1_x$.
Then there exists a unique density matrix $\rho_y$ such that
\be
\omega_y(A)&=&\Tr\rho_y A^{(1)},
\qquad A\in{\cal A}.
\nnee
It is shown in \cite{NJ19b} that the operator $X_{y,x}$ 
is given by
\be
X_{y,x}
&=&
S_x\left[\rho_x^{-1}\rho_y\otimes\Io\right]S_x\cr
&=&
J_x\left(\rho_x^{-1/2}\rho_y\rho_x^{-1/2}\otimes\Io\right)J_x
\nnee
and satisfies 
\be
X_{y,x}x&=&\left[\rho_y\rho_x^{-1}\otimes\Io\right]x.
\nnee
The vector $y$ is then given by $y=X_{y,x}^{1/2}x$.

\paragraph{Exponential arcs}

Take $y$ in ${\mathscr C}^1_x$. One has, making use of $J_x^2=\Io$ and $J_xx=x$,
\be
X^{t/2}_{y,x}x
&=&
J_x\left([\rho_x^{-1/2}\rho_y\rho_x^{-1/2}]^{t/2}\otimes\Io\right)J_xx\cr
&=&
\sum_i\sqrt{p_i}
J_x\left([\rho_x^{-1/2}\rho_y\rho_x^{-1/2}]^{t/2}e_i\otimes e_i\right).
\nnee
Hence one finds
\be
\zeta_{y,x}(t)
&=&\log ||X^{t/2}_{y,x}x||\cr
&=&
\log||\sum_i\sqrt{p_i}[\rho_x^{-1/2}\rho_y\rho_x^{-1/2}]^{t/2}e_i\otimes e_i||\cr
&=&
\frac 12\log\sum_i p_i
\left([\rho_x^{-1/2}\rho_y\rho_x^{-1/2}]^{t}e_i,e_i\right)\cr
&=&
\frac 12\log \Tr \rho_x [\rho_x^{-1/2}\rho_y\rho_x^{-1/2}]^{t}.
\nnee

Take $y$ in ${\mathscr C}^1_x$.
The exponential arc $t\mapsto \omega_t$ connecting $\omega_y$ to $\omega_x$ is given by
\be
\omega_t(A)
&=&
e^{-2\zeta_{y,x}(t)}(AX^{t/2}_{y,x}x,X^{t/2}_{y,x}x))\cr
&=&
e^{-2\zeta_{y,x}(t)}
\sum_{i,j}\sqrt{p_ip_j}
\bigg(AJ_x\left([\rho_x^{-1/2}\rho_y\rho_x^{-1/2}]^{t/2}e_i\otimes e_i\right),\cr
& &\quad
J_x\left([\rho_x^{-1/2}\rho_y\rho_x^{-1/2}]^{t/2}e_j\otimes e_j\right)\bigg)\cr
&=&
e^{-2\zeta_{y,x}(t)}
\sum_{i,j}\sqrt{p_ip_j}
\bigg(
[\rho_x^{-1/2}\rho_y\rho_x^{-1/2}]^{t/2}e_j\otimes e_j,\cr
& &\quad
A'\left([\rho_x^{-1/2}\rho_y\rho_x^{-1/2}]^{t/2}e_i\otimes e_i\right)
\bigg)
\nnee
with $A'=J_xAJ_x$. Use that $A'$ belongs to the commutant ${\cal A}'$ to write
$A'=\Io\otimes R$. One obtains
\be
\omega_t(A)
&=&
e^{-2\zeta_{y,x}(t)}
\sum_{i,j}\sqrt{p_ip_j}
\bigg(
[\rho_x^{-1/2}\rho_y\rho_x^{-1/2}]^{t/2}e_j\otimes e_j,\cr
& &\qquad\qquad
\left([\rho_x^{-1/2}\rho_y\rho_x^{-1/2}]^{t/2}e_i\otimes Re_i\right)
\bigg)\cr
&=&
e^{-2\zeta_{y,x}(t)}
\sum_{i,j}\sqrt{p_ip_j}
\left(
[\rho_x^{-1/2}\rho_y\rho_x^{-1/2}]^{t/2}e_j,
[\rho_x^{-1/2}\rho_y\rho_x^{-1/2}]^{t/2}e_i\right)\cr
& &
\times(e_j,Re_i).
\nnee
Now use that (see the appendix of \cite{NJ18})
\be
(e_j,Re_i)&=&(e_i,A^*e_j)
\nnee
to obtain 
\be
\omega_t(A)
&=&
e^{-2\zeta_{y,x}(t)}
\sum_{i,j}\sqrt{p_j}
\left(
[\rho_x^{-1/2}\rho_y\rho_x^{-1/2}]^{t/2}e_j,
[\rho_x^{-1/2}\rho_y\rho_x^{-1/2}]^{t/2}e_i\right)\cr
& &
\times(e_i,\rho_x^{1/2}A^*e_j)\cr
&=&
e^{-2\zeta_{y,x}(t)}
\sum_{j}
\left(
[\rho_x^{-1/2}\rho_y\rho_x^{-1/2}]^{t/2}\rho^{1/2}_xe_j,
[\rho_x^{-1/2}\rho_y\rho_x^{-1/2}]^{t/2}\rho_x^{1/2}A^*e_j\right)\cr
&=&
\Tr\rho_t A
\nnee
with the density matrix $\rho_t$ given by
\be
\rho_t
&=&
e^{-2\zeta_{y,x}(t)}\rho_x^{1/2}[\rho_x^{-1/2}\rho_y\rho_x^{-1/2}]^{t}\rho^{1/2}_x.
\label{dens:result}
\ee
This exponential arc differs from the geodesic studied in \cite{NJ18,NJ19b}.
The latter is given by
\be
\rho_t&\sim&\exp((1-t)\log\rho_x+t\log\rho_y).
\label{dens:altresult}
\ee
The problem when generalizing (\ref {dens:altresult})
to infinite dimensions is that the operators
$\log\rho_x$ and $\log\rho_y$ need not to commute and even need not to have a common domain
so that the sum appearing in the exponential function of (\ref {dens:altresult}) may be ill-defined.
Note that if $\rho_x$ and $\rho_y$ commute then (\ref {dens:result}) and (\ref {dens:altresult})
coincide.

\subsection{Classical probability}

\label{sect:abel}

In this section the von Neumann algebra $\cal A$ is (isomorphic with) the space
$L_\infty(\Ro^n,\Co)$ of all essentially bounded
complex functions on $\Ro^n$ with its Lebes\-gue measure. The Hilbert space $\cal H$
coincides with the space of square-integrable complex functions ${\mathscr L}_2(\Ro^n)$.
The symbol $x$ denotes an element of $\Ro^n$ instead of a vector in $\cal H$.
A function $A(x)$, element of $\cal A$,
acts on a square-integrable function $f(x)$ by pointwise multiplication $(Af)(x)=A(x)f(x)$.
The commutant ${\cal A}'$ of $\cal A$ coincides with $\cal A$.
See \cite {DJ81}, Part I, Ch. 7, Thm. 2.

Any normalized element $f$ of ${\mathscr L}_2(\Ro^n)$ determines a state $\omega_f$ on $\cal A$
by
\be
\omega_f(A)&=&\int_{\Ro^n}A(x)|f(x)|^2\upd x.
\nnee
The correspondence between probability measures $|f(x)|^2\upd x$
and vector states $\omega_f$ is one-to-one.

The function $f$ is cyclic and separating for $\cal A$ if it is strictly positive almost everywhere.
The cone ${\mathscr C}_f$ consists of all square integrable functions $g$ for which the ratio
$\tilde g=g/f$ is a non-negative function a.e.. The operator $X^{1/2}_{g,f}$ equals the function 
$\tilde g$. Clearly is $X^{1/2}_{g,f}f=g$. The state $\omega_g$ satisfies
$\omega_g(A)=\omega_f(X_{g,f}A)$. The operator $X_{g,f}$ is the Radon–Nikodym derivative
of the measure $|g(x)|^2\upd x$ w.r.t.~the measure $|f(x)|^2\upd x$. 

The square integrable function $g$ belongs to $\breve{\mathscr C}_f$ if $\tilde g$ is a bounded function.

\paragraph{Exponential arcs}

Take  $g$ in ${\mathscr C}^1_f$.
The exponential arc $\gamma_{g,f}$ connecting $g$ to $f$ is given by
\be
\gamma_{g,f}(t)(x)&=&e^{-\zeta_{g,f}(t)}[\tilde g(x)]^tf(x),
\nnee
with the normalization $\zeta_{g,f}$ given by
\be
\zeta_{g,f}(t)&=&\frac 12\log\int_{\Ro^n}[\tilde g(x)]^{2t}|f(x)|^2\upd x.
\nnee
The Fr\'echet derivative satisfies the equation
\be
\dot\gamma_{g,f}(t)(x)
&=&
\left[-\zeta'_{g,f}(t)+\log\tilde g(x)\right]\gamma_{g,f}(t)(x)
\nnee
with $\zeta'_{g,f}(t)$ the derivative of $\zeta_{g,f}(t)$.
It is well-defined when $\tilde g(x)=0$ implies $\gamma_{g,f}(t)(x)=0$.
In particular, the tangent vector at $t=0$ is given by
\be
\dot\gamma_{g,f}(0)(x)
&=&
\left[\log\tilde g(x)-\omega_f(\log\tilde g(x))\right]f(x).
\nnee

The exponential arc $\gamma_{g,f}(t)$ is discontinuous at $t=0$ when $\tilde g(x)$ is not strictly
positive a.e.. Hence, the set $\mathring\gamma_f$ consist of all vectors $g$ in ${\mathscr C}^1_f$
for which $\tilde g(x)>0$ a.e.. The corresponding probability measures 
$|f(x)|^2\upd x$ and $|g(x)|^2\upd x$ are then equivalent measures.
The set $\Vo_x$ contains all states $\omega_g$ such that 
the corresponding probability measures are absolutely continuous w.r.t.~the measure $|f(x)|^2\upd x$.

The exponential arc connecting the state $\omega_g$ to the state $\omega_f$ is given by
\be
\omega_t(A)
&=&\int_{\Ro^n}A(x)|\tilde g(x)|^{2t}\,|f(x)|^2\upd x\cr
&=&
\int_{\Ro^n}A(x)\,|f(x)|^{2(1-t)}\,|g(x)|^{2t}\upd x.
\nnee
This expression for the exponential arc is found for instance in \cite{PG13}.

\section{Discussion}

The main result of the present paper is the definition of exponential arcs
both in Hilbert space $\mathscr H$  and in the space $\Vo$
of vector states on a $\sigma$-finite von Neumann algebra $\cal A$ in its standard form \cite{HU75}.
A proof is given  that these definitions are logically consistent.
See Theorems  \ref {exparchil:prop:subarc} and \ref {exparcstate:prop:subarc}.
Theorems \ref {tangent:thm} and \ref {neigh:thm} show that the arcs are Fr\'echet differentiable.
The special cases treated in the previous section show that the formulation in terms
of $\sigma$-additive von Neumann algebras covers both quantum models and models of classical
probability.

The present work makes use of 
powerful results of Tomita-Takesaki theory \cite{TM70}.
The theory describes the symmetry that may exist between a von Neumann algebra and its commutant.
The underlying structure is a symmetry between a pair of real subspaces of the Hilbert space \cite{RvD77}.
This raises the hope that eventually one can formulate a theory in which unbounded operators
can be avoided.

The generalization of the notion of an exponential arc to a non-commutative context
can be done in more than one way.
The exponential arcs introduced in the present paper differ from the geodesics
discussed in \cite{NJ18,NJ19b}.
Araki \cite{AH74} introduces an $\alpha$-family of cones in Hilbert space, with $\alpha$
running from $0$ to $1/2$. The present paper focuses on the $\alpha=1/2$-cone while
\cite{NJ18,NJ19b} is intended to correspond with the self-dual case $\alpha=1/4$. 
Some of the technical problems
one encounters in the latter case, when trying to generalize to an infinite-dimensional
Hilbert space, are avoided in the present paper.
In the commutative case all $\alpha$-cones coincide and the present approach reproduces
known results.

The definition of exponential arcs is only a first step in the study of the manifold of
vector states from the view point of Information Geometry. The next step is the study
of tangent planes, the search for Banach or Hilbert charts and the introduction of a
Riemannian metric.

\section*{Acknowledgment}
I am indebted to an anonymous referee carefully reading a previous version of the paper
and suggesting to work out the link with the natural positive cone.
This lead to a complete overhaul of the paper.

\end{document}